\documentclass{article}

\usepackage{amsfonts,graphicx,amsmath,subfigure}

\def\R{\mathbb R}

\def\epsilon{\varepsilon}
\def\ds{\displaystyle}

\newcommand{\be}{\par\nobreak\noindent  \begin{equation}}
\newcommand{\ee}{\end{equation}}
\newcommand{\baa}{\begin{array}}
\newcommand{\eaa}{\end{array}}
\newcommand{\ba}{\begin{eqnarray}}
\newcommand{\ea}{\end{eqnarray}}

\newtheorem{lemma}{Lemma}[section]
\newtheorem{theorem}[lemma]{Theorem}
\newtheorem{corollary}[lemma]{Corollary}

\newtheorem{proposition}[lemma]{Proposition}
\newtheorem{remark}[lemma]{Remark}

\title{Biological invasions: deriving the regions at risk from partial measurements}

\author{Michel Cristofol$^{\hbox{\small{a}}}$ and Lionel Roques$^{\hbox{\small{ b,* }}}$ \\
\footnotesize{$^{\hbox{a} }$ LATP, CNRS/UMR 6632, CMI,
Universit\'e de Provence, France} \\
\footnotesize{and Universit\'e d'Aix-Marseille III, IUT de St
J\'er\^ome, France} \\
\footnotesize{$^{\hbox{b }}$ INRA, UR546 Biostatistique et Processus Spatiaux, F-84914 Avignon, France} \\
\footnotesize{$^{\hbox{*} }$ Corresponding author. E-mail:
lionel.roques@avignon.inra.fr} }

\date{}

\begin{document}

\maketitle

\begin{abstract}
We consider the problem of forecasting the regions at higher risk
for  newly introduced invasive species. Favourable and
unfavourable regions may indeed not be known a priori, especially
for exotic species whose hosts in native range and newly-colonised
areas can be different. Assuming that the species is modelled by a
logistic-like reaction-diffusion equation, we prove that the
spatial arrangement of the favourable and unfavourable regions can
theoretically be determined using only
partial measurements of the population density: 1) a local
``spatio-temporal" measurement, during a short time period and, 2)
a ``spatial" measurement in the whole region susceptible to
colonisation. We then present a stochastic algorithm which is proved
analytically, and then on several numerical examples, to be
effective in deriving these regions.
\end{abstract}


\section{Introduction}\label{section1}

Because of trade globalisation,  a substantial increase in
biological invasions has been observed over the last decades (e.g.
Liebhold et al. \cite{lieb}). These invasive species are, by
definition \cite{nisic}, likely to cause economic or environmental
harm or harm to human health. Thus, it is a major concern to
forecast, at the beginning of an invasion, the areas which will be
more or less infested by the species.

Because of their exotic nature, invading species generally face
little competition or predation. They are therefore well adapted
to modelling via single-species models.

Reaction-diffusion models have proved themselves to give good
qualitative results regarding biological invasions (see the
pioneering paper of Skellam \cite{ske}, and the books \cite{sk},
\cite{tur} and \cite{oku} for review).

The most widely used single-species reaction-diffusion model, in
homogeneous environments, is probably the Fisher-Kolmogorov \cite{fi,kpp} model:
\begin{equation}\label{kpp}
u_t=D \Delta u+u(\mu-  \gamma u), \ t>0,  \ x\in \Omega \subset
\R^N ,
\end{equation}
where $u=u(t,x)$ is the population density at time $t$ and space
position $x$, $D$ is the diffusion coefficient,   $\mu$
corresponds to the \emph{constant} intrinsic growth rate, and
$\frac{\mu}{\gamma}$ is the environment's carrying capacity. Thus
$\gamma$ measures the susceptibility to crowding effects.

On the other hand, the environment is generally far from being
homogeneous. The spreading speed of the invasion, as well as the
final equilibrium attained by the population are in fact often
highly dependent on these heterogeneities (\cite{sk}, \cite{bhr1},
\cite{bhr2}, \cite{kks}). A natural extension of (\ref{kpp}) to
heterogeneous environments has been introduced by Shigesada,
Kawasaki, Teramoto \cite{skt}: \be\label{skt} u_t=
\nabla(D(x)\nabla u)+u(\mu(x)-\gamma(x) u), \ t>0,  \ x\in \Omega
\subset \R^N.\ee In this case, the diffusivity matrix $D(x)$, and
the coefficients $\mu(x)$ and $\gamma(x)$  depend on the space
variable $x$, and can therefore include some effects of
environmental heterogeneity.

In this paper, we consider the simpler case where $D(x)$  is
assumed to be constant and isotropic and $\gamma$ is also assumed to be
positive and constant: \be\label{mod1} u_t= D\Delta
u+u(\mu(x)-\gamma u), \ t>0,  \ x\in \Omega \subset \R^N.\ee The
regions where $\mu$ is high correspond to favourable regions (high
intrinsic growth rate and high environment carrying capacity),
whereas the regions with low values of $\mu$ are less favourable,
or even unfavourable when $\mu<0$. In what follows, in order to
obtain clearer biological interpretations of our results, we say
that $\mu$ is a ``habitat configuration".

With this type of model, many qualitative results have been
established, especially regarding the influence of spatial
heterogeneities of the environment on population persistence, and
on the value of the equilibrium population density (\cite{sk},
\cite{bhr1}, \cite{ccL},   \cite{rs}, \cite{rc1}). However, for a newly
introduced species, like an invasive species at the beginning of
its introduction, the regions where $\mu$ is high or low may not
be known \textit{a priori}, particularly when the environment is
very different from that of the species native range.

In this paper, we propose a method of deriving the habitat
configuration $\mu$, basing ourselves only on partial measurements
of the population density at the beginning of the invasion
process. In section \ref{section2}, we begin by giving a precise
mathematical formulation of our estimation problem. We then
describe our main mathematical results, and we link them with
ecological interpretations. These theoretical results form the
basis of an algorithm that we propose, in section \ref{section3},
for recovering the habitat configuration $\mu$. In section
\ref{section4}, we provide numerical examples illustrating our
results. These results are further discussed in section
\ref{section5}.

\section{Formulation of the problem and main results}\label{section2}

\subsection{Model and hypotheses}\label{section2.1}

We assume that the population density $u_\gamma$ is governed by
the following parabolic equation:
\begin{equation*} \left\{ \baa{l}  \partial_t
u_\gamma = D \Delta u_\gamma + u_\gamma (\mu(x)-\gamma u_\gamma), \;\; t>0, x \in \Omega,\\
u_\gamma(t,x)=0, \;\; t>0, x \in \partial \Omega, \\ u_\gamma(0,x)=u_{i}(x)
\hbox{ in }\Omega, \eaa \right. \ (P_{\mu,\gamma})\end{equation*}
where $\Omega$ is a bounded subdomain of $\R^d$ with
boundary $\partial \Omega$. We will denote $Q:=(0,+\infty)\times \Omega$ and
$\Sigma:=(0,+\infty)\times \partial \Omega$.

\

The growth rate function $\mu$ is \textit{a priori} assumed to be
bounded, and to take a known constant value outside a fixed compact subset $\Omega_1$
of $\Omega$:
$$\mu \in \mathcal{M}:=\{\rho\in  L^{\infty}(\Omega), \ -M\leq\rho\leq
M \hbox{ a.e., and }\rho\equiv m \hbox{ in } \Omega\backslash
\Omega_1\},$$for some  constants $m,M$, with $M>0$; the notation ``a.e.'' means  ``almost everywhere'', which is equivalent to ``except on a set of zero measure''.

The initial population density $u_{i}(x)$ is assumed to be bounded
(in $C^2(\overline{\Omega})$), and bounded from below by a fixed
positive constant in a fixed closed ball $\mathcal{B}_{\epsilon} \subset
 \Omega_1 $, of small radius $\epsilon$: \be
\mathcal{D}:=\{\phi \ge 0, \; \phi \in C^2(\overline{\Omega}),  \ \| \phi
\|_{C^2(\Omega)} \leq \overline{u_{i}}, \ \phi \geq
\underline{u_{i}} \mbox{ in } \mathcal{B}_{\epsilon}, \}, \ee for
some positive constants $\overline{u_{i}}$ and
$\underline{u_{i}}$.

Absorbing (Dirichlet)  boundary conditions are assumed.

\begin{remark}
{\rm
Absorbing boundary conditions mean that the individuals crossing the
boundary immediately die. Such conditions can be ecologically relevant in numerous situations. For instance for many plant species, seacoasts are lethal and thus constitute this kind of boundaries.

For technical reasons we have to introduce the subset $\Omega_1$, such that, in the interface between $\Omega_{1}$ and $\Omega$, $\mu$ takes a known value $m$. This value is typically negative, indicating that, near the lethal boundary, the environment is unfavourable.  This assumption is not very restrictive since, in fact, $\Omega_1$ can be chosen as close as we want to $\Omega$.

For  precise definitions of the functional spaces $L^{2}$, $L^{\infty}$ and $C^2$ as well as the other mathematical notations used throughout this paper, the reader can refer, e.g., to \cite{evans}.}
\end{remark}

\subsection{Main question}\label{section2.2}

The main question that we presented at the end of the Introduction
section can now be stated: for any time-span $(t_0,t_1)$, and any non-empty
subset $\omega$ of $\Omega_1$,  is it possible to
estimate the function $\mu(x)$ in $\Omega$, basing ourselves only
on measurements of $u_{\gamma}(t,x)$ over $ (t_0,t_1)\times \omega$, and on
a single measurement of $u_{\gamma}(t,x)$ in the whole domain $\Omega$ at a
time $T'=\frac{t_0+t_1}{2}$?

\subsection{Estimating the habitat configuration}\label{section2.3}

Let $\tilde{\mu}$ be a function in $\mathcal{M}$, and  let
$\tilde{v}$ be the solution of the linear parabolic problem
$(P_{\tilde{\mu},0})$.  We define a  functional $G_{\mu}$, over
$\R_+\times \mathcal{M}$, by \ba G_{\mu}(\gamma,\tilde{\mu}) & = &
\|\partial_t
u_{\gamma}-\partial_t\tilde{v}\|^2_{L^2((t_0,t_1)\times\omega)}\\
& &+ \|\Delta u_{\gamma}\left( T' , \cdot\right)  -\Delta
\tilde{v}\left( T' , \cdot\right)\|^2_{L^2(\Omega)}+ \|
u_{\gamma}\left( T' , \cdot\right)- \tilde{v}\left( T' ,
\cdot\right)\|^2_{L^2(\Omega)}, \ea where $u_{\gamma}$ is the
solution of $(P_{\mu,\gamma})$.
This functional $G_{\mu}$ quantifies the gap between $u_{\gamma}$ and $\tilde{v}$ on the set where $u_{\gamma}$ has been measured.

\begin{theorem}
The functions $\mu,\tilde{\mu}\in \mathcal{M}$ being given, we
have:
$$\|\mu-\tilde{\mu}\|^2_{L^2(\Omega_1)}\leq \frac{C}{\overline{u_{i}}^2} G_{\mu}(0,\tilde{\mu}),$$for all $\tilde{\mu}\in \mathcal{M}$
and for some positive constant
$C=C(\Omega,\Omega_1,\omega,\mathcal{B}_{\epsilon},D,t_0,t_1,\underline{u_i}/\overline{u_i})$.
\label{th1}
\end{theorem}
The proof of this result is given in Appendix A.  It bears on a
Carleman-type estimate.

\textbf{Biological interpretation:}  This stability result means
that, in the linear case corresponding to Malthusian populations
($\gamma=0$), two different habitat configurations $\mu, \
\tilde{\mu}$ cannot lead to close population densities $u_0, \
\tilde{v}$. Indeed, having
population densities that are close to each other in the two
situations,  even on a very small region $\omega$, during a small
time period $(t_0,t_1)$, and in the whole space $\Omega$ at a
single time $T'$, would lead to small $G_\mu$ values, and
therefore, from Theorem \ref{th1}, to close values of the growth
rate coefficients $\mu$ and $\tilde{\mu}$.

Theorem \ref{th1} implies the following uniqueness result:
\begin{corollary}
If $v$ is a solution of both $(P_{\mu,0})$ and
$(P_{\tilde{\mu},0})$, then $\mu=\tilde{\mu} $ a.e.  in $\Omega_1$, and
therefore in $\Omega$. \label{cor1}
\end{corollary}

\textbf{Biological interpretation:}  In the linear case
($\gamma=0$), if two habitat configurations $\mu, \ \tilde{\mu}$
lead to identical
population densities $u_0$, $\tilde{v}$,   even on a very small
region $\omega$, during a small time period $(t_0,t_1)$, and in
the whole space $\Omega$ at a single time $T'$, then these habitat
configurations are identical.

Next we have the following result:
\begin{theorem}\label{th2}
We have\footnote{Two functions
$f(\mu,\tilde{\mu},u_i,\underline{u_i},\overline{u_i},\gamma)$ and
$g(\mu,\tilde{\mu},u_i,\underline{u_i},\overline{u_i},\gamma)$,
are written $f=\mathcal{O}(g)$ as $g\to 0$ if there exists a
constant $K>0$, independent of $\mu$, $\tilde{\mu}$, $u_i$,
$\underline{u_i},$ $\overline{u_i}$ and $\gamma$, such that
$|f|\leq K |g|$ for $g$ small enough.}
$|G_{\mu}(0,\tilde{\mu})-G_{\mu}(\gamma,\tilde{\mu})|=\mathcal{O}(\overline{u_{i}}^3),$
as $\overline{u_i}\to 0.$
\end{theorem}
The proof of this result is given in Appendix B.

\textbf{Biological interpretation:} Assume that the habitat
configuration $\mu$ is not known, but that we have measurements of
the population density $u_\gamma$, governed by the full nonlinear
model (\ref{mod1}).  Consider a configuration
$\tilde{\mu}$ in $\mathcal{M}$
 such that the population density $\tilde{v}$ obtained as a solution of the linear model
$(P_{0,\tilde{\mu}})$ has values close to those taken by the
population density $u_\gamma$, in the sense that
$G_{\mu}(\gamma,\tilde{\mu})$ is close to $0$.
If the initial population density is far from
the environment carrying capacity, then
$\overline{u_i} \ll \frac{\mu}{\gamma}$, $\overline{u_{i}}$ is small
and, from Theorem \ref{th2}, $G_{\mu}(0,\tilde{\mu})$ is also
close to $0$. Thus Theorem \ref{th1} implies that the habitat
configuration $\tilde{\mu}$ is an accurate estimate of $\mu$. In
section \ref{section3}, we propose an algorithm to obtain
explicitly such estimates of $\mu$.

\begin{remark}
{\rm In fact, the term $\mathcal{O}(\overline{u_{i}}^3)$ increases exponentially with time $t_{1}$. Thus, obtaining accurate estimates of $\mu$ require, in practice, to work with small times i.e. at the beginning of the invasion.}
\end{remark}

\subsection{Forecasting the fate of the invading population}\label{section2.4}

The knowledge of an $L^2$-estimate $\tilde{\mu}$ of $\mu$
enables us to give an estimate of the asymptotic behaviour of the
solution $u_{\gamma}$ of $(P_{\mu,\gamma})$, as $t \to +\infty$,
and especially to know whether the population will become extinct or not.
Indeed, as $t\to +\infty$, it is known that (see e.g. \cite{bhr1},
for a proof with another type of boundary condition) the solution
$u_\gamma(t,x)$ of $(P_{\mu,\gamma})$ converges to the unique
nonnegative and bounded solution $p_\gamma$ of
\begin{equation*} -D \Delta p_\gamma = p_\gamma (\mu(x)-\gamma p_\gamma) \hbox{ in }\Omega, \ \ \ \ (S_{\mu,\gamma})\end{equation*}
with $p_\gamma=0$ on $\partial \Omega$. Moreover, $p_\gamma\equiv
0$ if and only if $\lambda_1[\mu]\geq 0$, where $\lambda_1[\mu]$
is the smallest eigenvalue of the elliptic operator $\mathcal{L}:
\ \psi \mapsto - D \Delta \psi -\mu(x) \psi$, with Dirichlet
boundary conditions. On the other hand, if $\lambda_1[\mu]<0$,
then $p_\gamma \equiv
0$ in $\Omega$ (note that $\gamma$ does not appear in the
definition of $\lambda_1$).

We have the following result.
\begin{proposition}
Let us consider a sequence $(\tilde{\mu}_n)_{n\geq 0}$ in
$\mathcal{M}$, such that $\tilde{\mu}_n\to \mu$ in $L^2(\Omega)$
as $n\to \infty$.

a) The solution $\tilde{p}_{\gamma,n}$ of the problem
$(S_{\tilde{\mu}_n,\gamma})$ converges to $p_\gamma$ as $n\to
\infty$, uniformly in $\Omega$.

b) $\lambda_1[\tilde{\mu}_n]\to \lambda_1[\mu]$ as $n\to +\infty$.

\label{lemp}
\end{proposition}
The proof of this result is classical and can be found in
\cite{bhr1,rh}.

\textbf{Biological interpretation:} Assume that the habitat
configuration $\mu$ is not known. We know that, in large times,
the population density $u_\gamma$ will tend to an unknown steady
state $p_\gamma$ (possibly $0$, in case of extinction of the
population). The part a) of the above proposition means that, if
we know an accurate ($L^2$-) estimate $\tilde{\mu}$ of $\mu$, then
we can deduce an accurate estimate $\tilde{p}_\gamma$ of the
steady state $p_\gamma$, provided the coefficient $\gamma$ is
known. Part b) shows that, even if $\gamma$ is not known, having
an estimate $\tilde{\mu}$ of $\mu$ enables to obtain an estimate
of $\lambda_1[\mu]$, and therefore to forecast whether the species
will survive or not. Indeed the sign of $\lambda_1[\mu]$ controls
the fate of the invading species (persistence if
$\lambda_1[\mu]<0$ and extinction if $\lambda_1[\mu]\ge 0$, see
\cite{bhr1,ccL,rs,rh} for more details) .

\section{Simulated annealing algorithm}\label{section3}

Let  $(t_0,t_1)$ be a fixed time interval, and $\omega
\subset \Omega_1$ be fixed. We assume that we have measurements of
the solution $u_{\gamma}(t,x)$ of $(P_{\mu,\gamma})$  over $
(t_0,t_1)\times \omega$, and of
$u_{\gamma}(\frac{t_0+t_1}{2},x)$ in $\Omega$.
However, the function $\mu$ and the constant $\gamma$ are assumed to
be unknown. Our objective is to build an algorithm for recovering
$\mu$.
\begin{remark}
{\rm When the function $u_\gamma$ is known, the computation of
$G_{\mu}(\gamma,\cdot)$ does not require the knowledge of
$\gamma$.}
\end{remark}
The function $\mu$ is assumed to belong to a known finite subset
$E$ of $\mathcal{M}$, equipped with a neighbourhood system. We
build a sequence $\hat{\mu}_{n}$ of $N$ elements of $E$ with the
following simulated annealing algorithm:

\

\noindent $n=0$

\noindent Initialise $\hat{\mu}_{0}$

\noindent {\bf while} $n\leq N$

\noindent \ \ \ \ Choose randomly a neighbour $\nu$ of
$\hat{\mu}_{\gamma,n}$

\noindent \ \ \ \ {\bf if} $G_{\mu}(\gamma,\nu)\leq
G_{\mu}(\gamma,\hat{\mu}_{n})$

\noindent \ \ \ \ \ \ \ \ $\hat{\mu}_{n+1}\leftarrow\nu$

\noindent \ \ \ \ {\bf else}

\noindent \ \ \ \ \ \ \ \ Choose randomly with an uniform law
$w\in (0,1)$

\noindent \ \ \ \ \ \ \ \ {\bf if
}$w<e^{\frac{G_{\mu}(\gamma,\hat{\mu}_{n})-G_{\mu}(\gamma,\nu)}{\Theta(n)}}$

\noindent \ \ \ \ \ \ \ \ \ \ \ \ $\hat{\mu}_{n+1} \leftarrow\nu$

\noindent \ \ \ \ \ \ \ \ {\bf else }

\noindent \ \ \ \ \ \ \ \ \ \ \ \  $\hat{\mu}_{n+1}
\leftarrow\hat{\mu}_{n}$

\noindent \ \ \ \ \ \ \ \ {\bf endif }

\noindent \ \ \ \ {\bf endif}

\noindent \ \ \ \ $n \leftarrow n+1$

\noindent {\bf endwhile}

\

The sequence  $\Theta(n)$ (cooling schedule) is composed of real
positive numbers, decreasing to $0$. The simulated annealing
algorithm gives a sequence $\hat{\mu}_{n}$ of elements of $E$. It
is known (see e.g. \cite{hajek}) that, for a cooling schedule
$\Theta(n)$ which converges sufficiently slowly to $0$, this
sequence converges in $L^2(\Omega)$ to a global minimiser
$\hat{\mu}$ of $G_{\mu}(\gamma,.)$ in $E$ (but see Remark \ref{remcool}).

Moreover, from Theorems \ref{th1} and \ref{th2}, we have
$$\|\mu-\hat{\mu}_{n}\|^2_{L^2(\Omega)}\leq \frac{C}{\overline{u_{i}}^2} G_{\mu}(0,\hat{\mu}_{n})\leq  \frac{C}{\overline{u_{i}}^2} G_{\mu}(\gamma,\hat{\mu}_{n})+ \varepsilon_1(\overline{u_i}),$$
where $\varepsilon_1$ is a real-valued function such that
$\varepsilon_1(s)\to 0$ as $s\to 0$. Since $\mu\in E$ we obtain
that, for $n$ large enough,
$$G_{\mu}(\gamma,\hat{\mu}_{n})\leq G_{\mu}(\gamma,\mu),$$
and, from Appendix A
%(\ref{estimw3}),
$$  \frac{C}{\overline{u_{i}}^2} G_{\mu}(\gamma,\mu) \to 0 \hbox{ as }
\overline{u_i}\to 0,$$the ratio $\underline{u_{i}}/\overline{u_{i}}$ being kept constant. We finally get:
$$\|\mu-\hat{\mu}_{n}\|^2_{L^2(\Omega)}\to 0 \hbox{ as }
\overline{u_i}\to 0 \hbox{ and }n\to +\infty,$$for a fixed ratio $\underline{u_{i}}/\overline{u_{i}}$. Thus, for
$\overline{u_i}$ small enough, and for $n$ large enough,
$\hat{\mu}_{n}$ is as close as we want to $\mu$, in the
$L^2$-sense.

\begin{remark}
\label{remcool}{\rm
The cooling rate $\Theta(n)$ leading to the exact optimal configuration with probability $1$ decreases very slowly (logarithmically) and cannot be used in practice; see \cite{HJJ} for a detailed discussion. Empirically, a good trade-off between quality of solutions and time required for computation is obtained with exponential cooling schedules of the type $\Theta(n)=\Theta_0 \times \alpha^{n}$, with $\alpha<1$, first proposed by Kirkpatrick et al. \cite{Kirk}. Many other cooling schedules are possible, but too rapid cooling results in a system frozen into a state far from the optimal one. The starting temperature $\Theta_0$ should be chosen high enough to initially accept all changes $\hat{\mu}_{n+1} \leftarrow\nu$, whatever the  neighbour $\nu$.

For this type of algorithm, there are no general rules for the choice of the stopping criterion (see \cite{HJJ}), which should be heuristically adapted to the considered optimisation problem.}
\end{remark}

\section{Numerical computations}\label{section4}

In this section, in one-dimensional and two-dimensional cases, we
check that the algorithm presented in section \ref{section3} can
work in practice.

In each of the four following examples, we fixed the sets
$\Omega$, $\Omega_1$ and $\mathcal{M}$ and we defined a finite
subset $E\subset \mathcal{M}$ equipped with a neighbourhood system.
Then, for a fixed habitat configuration $\mu\in E$ we computed,
using a second-order finite elements method, the solution $u(t,x)$
of $(P_{\mu,\gamma})$, for $D=1$, $\gamma=0.1$, $t\in (0,0.5)$,
and for a given initial population density $u_i$. Then, we fixed
$t_0=0.1$, $t_1=0.4$, and a  subset $\omega \subset
\Omega_1$, and we stored the values of $u(\frac{t_0+t_1}{2},x)$,
for $x\in \Omega$ and $u(t,x)$ for $(t,x)\in (t_0,t_1)\times
\omega$. Using only these values, we computed the sequence
($\hat{\mu}_n$) of elements of $E$, defined by the simulated
annealing algorithm
of section \ref{section3}, the function $\hat{\mu}_0$
being sampled arbitrarily, in a uniform law, among the elements
of $E$.

In the following examples, we used $\Theta_{0}=100,$
and $\Theta(n)=100 \times 0.99^{n}$ for the cooling schedule.
Our stopping criterion was to have no change in the configuration $\hat{\mu}_{n}$ during $500$ iterations.

The rigourous definitions of the sets $E$ and of the associated neighbourhood systems which are used in the following examples can be found in Appendix C.

\subsection{One-dimensional case}

Assume that $\Omega=(0,100)$, $\Omega_1=[10,90]$,  $\omega=[55,58]$,   $M=2$,
$m=-1$ and $u_i=0.1 (1-x/100)
\sin(\pi x /25)^2$.

{\bf Example 1:} the set $E$ is composed of binary step functions taking only the values $m$ and $M$.

\begin{figure}
\centering
\subfigure[]{%
\includegraphics[width=6cm]{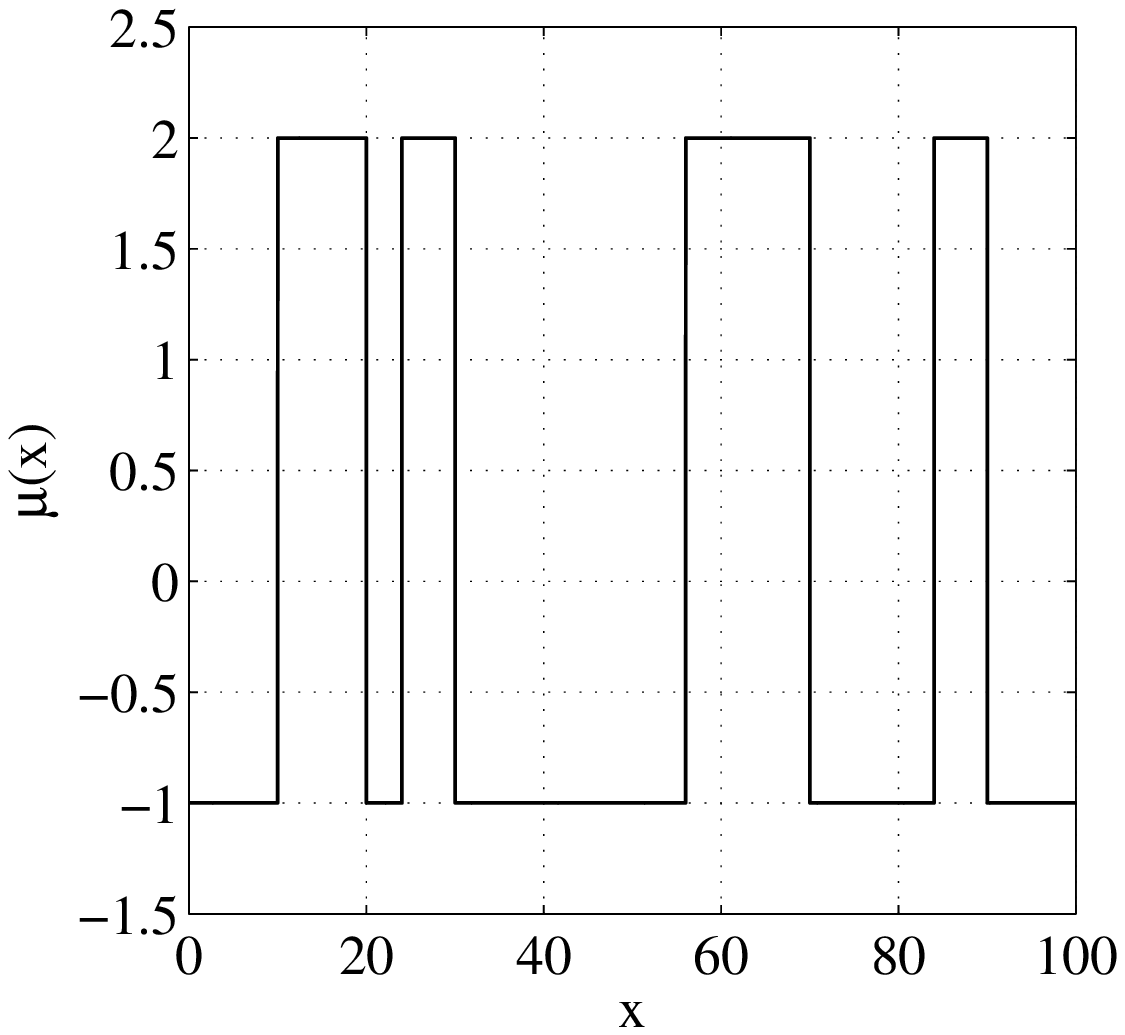}}
\subfigure[]{%
\includegraphics[width=6cm]{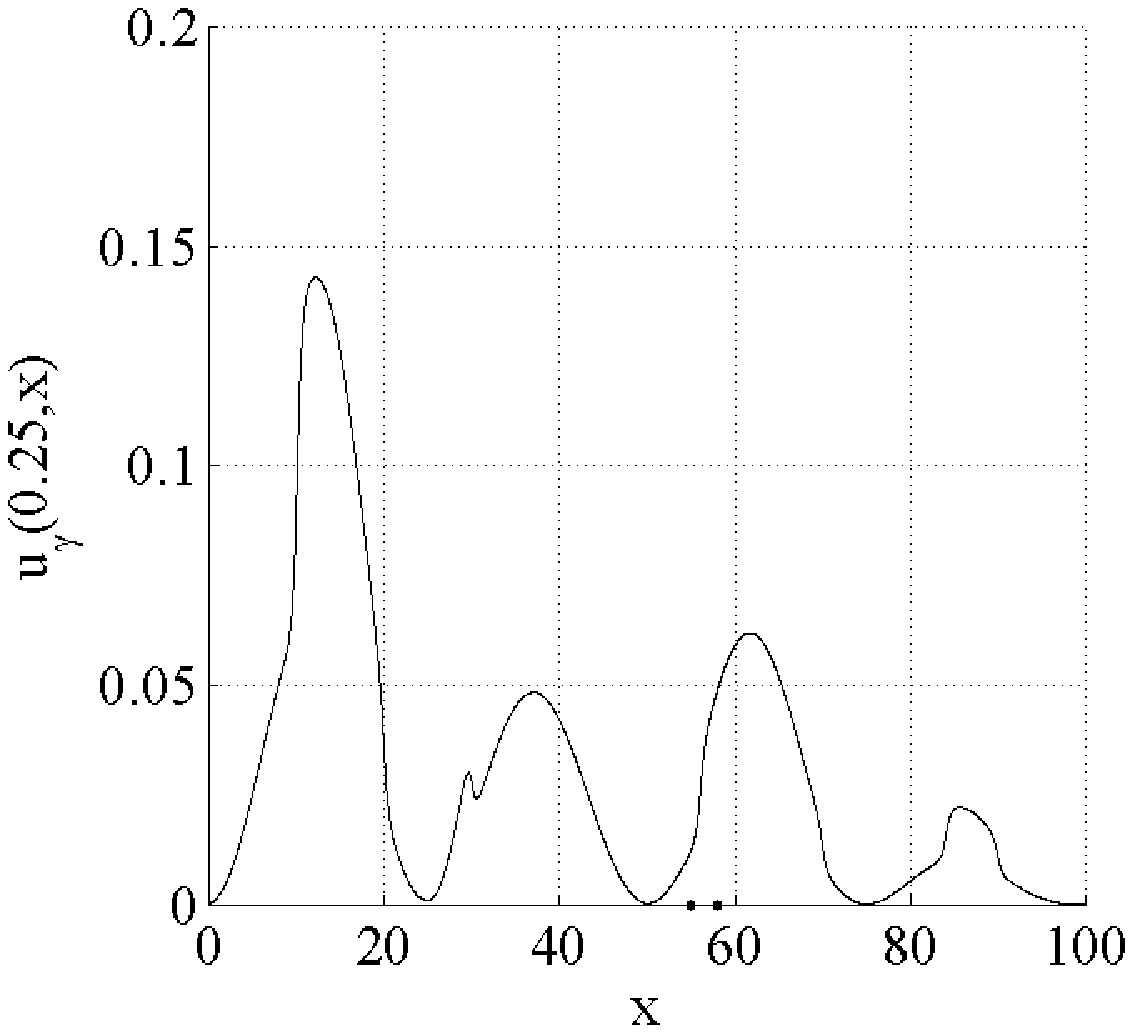}}
%\subfigure[]{%
%\includegraphics*[width=6cm]{exple1_mu700r.eps}}
%\subfigure[]{%
%\includegraphics*[width=6cm]{exple1_mu.eps}}
\caption{Example 1: a) exact habitat configuration $\mu$; b)
measurement of $ u_{\gamma}(0.25,x)$ which was used for recovering
$\mu$; the domain $\omega$ is delimited by two black dots. The
exact configuration $\mu$ was recovered after $1500$ iterations of
the algorithm.
}%
\label{fig:exple1}
\end{figure}

The function $\mu\in E$ and the measurement $ u_{\gamma}(\frac{t_0+t_1}{2},x)$ are depicted in Fig. \ref{fig:exple1}, (a) and (b).
Two elements of $E$ are said to be neighbours if they differ only on an interval of length $1$.

The
sequence ($\hat{\mu}_n$) stabilised  on the exact configuration
$\mu$ after about $N=1500$ iterations.

{\bf Example 2:} the set $E$ is composed of step functions which  can take 21 different values between $m$ and $M$.

\begin{figure}
\centering
\subfigure[]{%
\includegraphics[width=6cm]{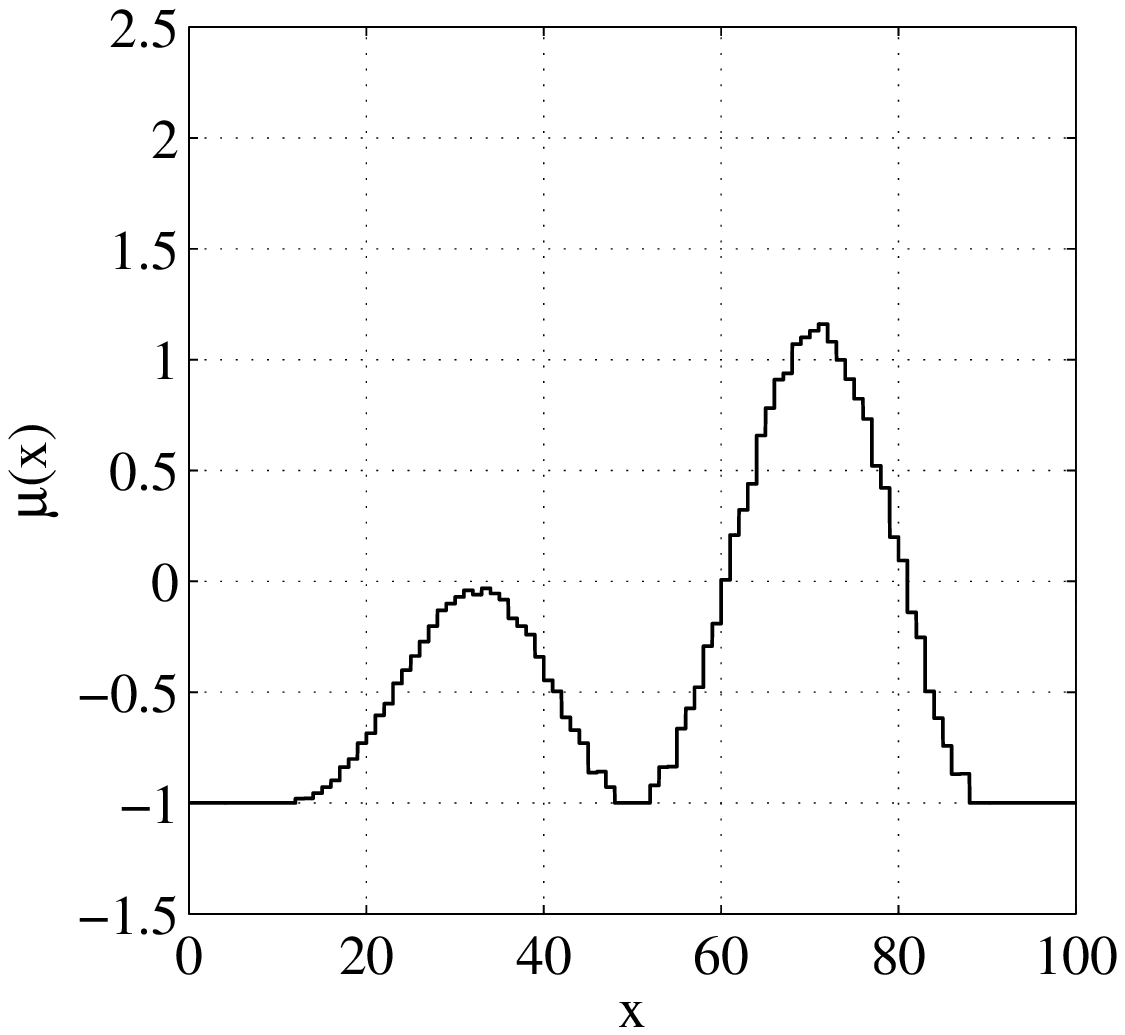}}
\subfigure[]{%
\includegraphics[width=6cm]{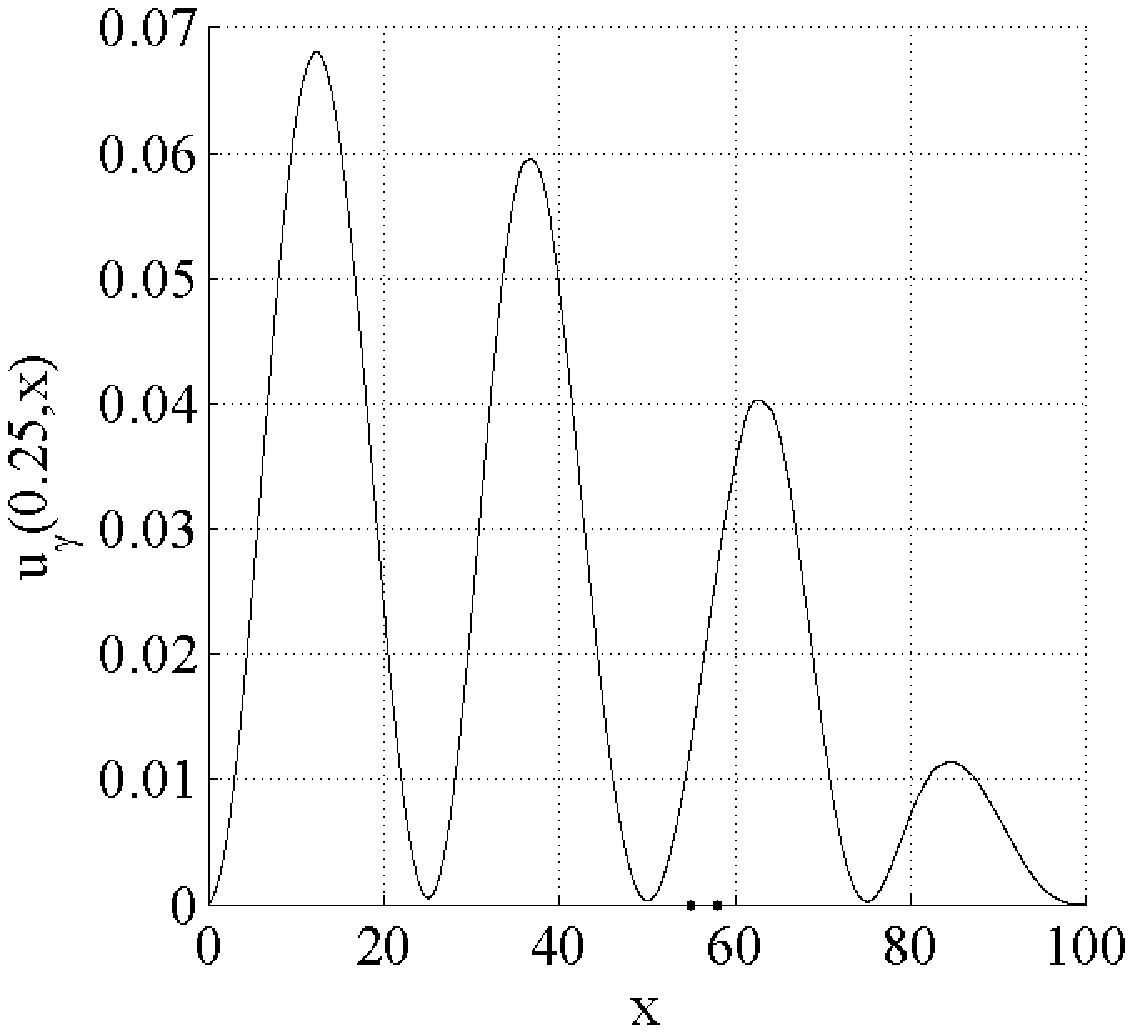}}
\subfigure[]{%
%\includegraphics*[width=6cm]{exple2_mu4000r.eps}}
%\subfigure[]{%
\includegraphics[width=6cm]{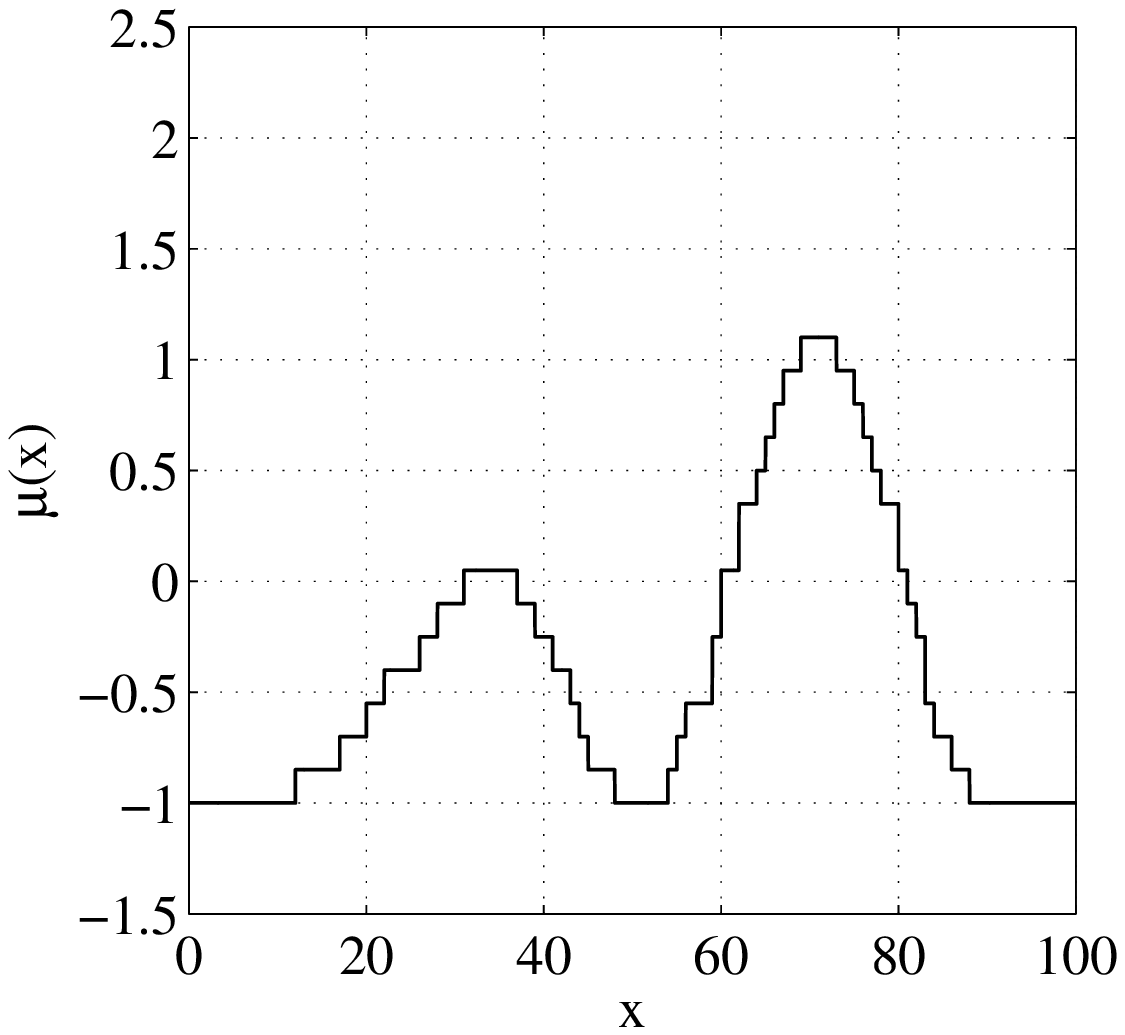}}
\caption{Example 2: a) exact habitat configuration $\mu$; b)
measurement of $ u_{\gamma}(0.25,x)$ which was used for recovering
$\mu$; the domain $\omega$ is delimited by two black dots; c)
configuration $\hat{\mu}:=\hat{\mu}_{7500}$,
obtained after $7500$ iterations.}%
\label{fig:exple2}
\end{figure}

The function $\mu\in E$ and the measurement $ u_{\gamma}(\frac{t_0+t_1}{2},x)$ are depicted in Fig. \ref{fig:exple2}, (a) and (b).
Two elements of $E$ are said to be neighbours if they differ by $(M-m)/20$ on an interval of length $1$.

This time, the sequence
($\hat{\mu}_n$) stabilised on a configuration $\hat{\mu}$ (Fig.
\ref{fig:exple2}, (c)) after $7500$ iterations. The mean error in
this case was
$\frac{1}{100}\int_0^{100}|\mu(x)-\hat{\mu}(x)|dx=0.05.$

\subsection{Two-dimensional case}

Assume now that $\Omega=(0,20)\times(0,20)$,
$\Omega_1=[2,18]\times[2,18]$, and that $\omega$
is the closed ball of centre $(7;7)$ and radius $3$ . Assume also that $M=2$ and $m=-1$, and that the initial data is $u_i=0.1 x y/400 \sin(x/4)^2 \sin(y/4)^2$.

{\bf Example 3:}  $E$ is composed of binary functions which can only take the values $m$ and $M$ on each cell of a regular lattice.

\begin{figure}
\centering
\subfigure[]{%
\includegraphics[width=6cm]{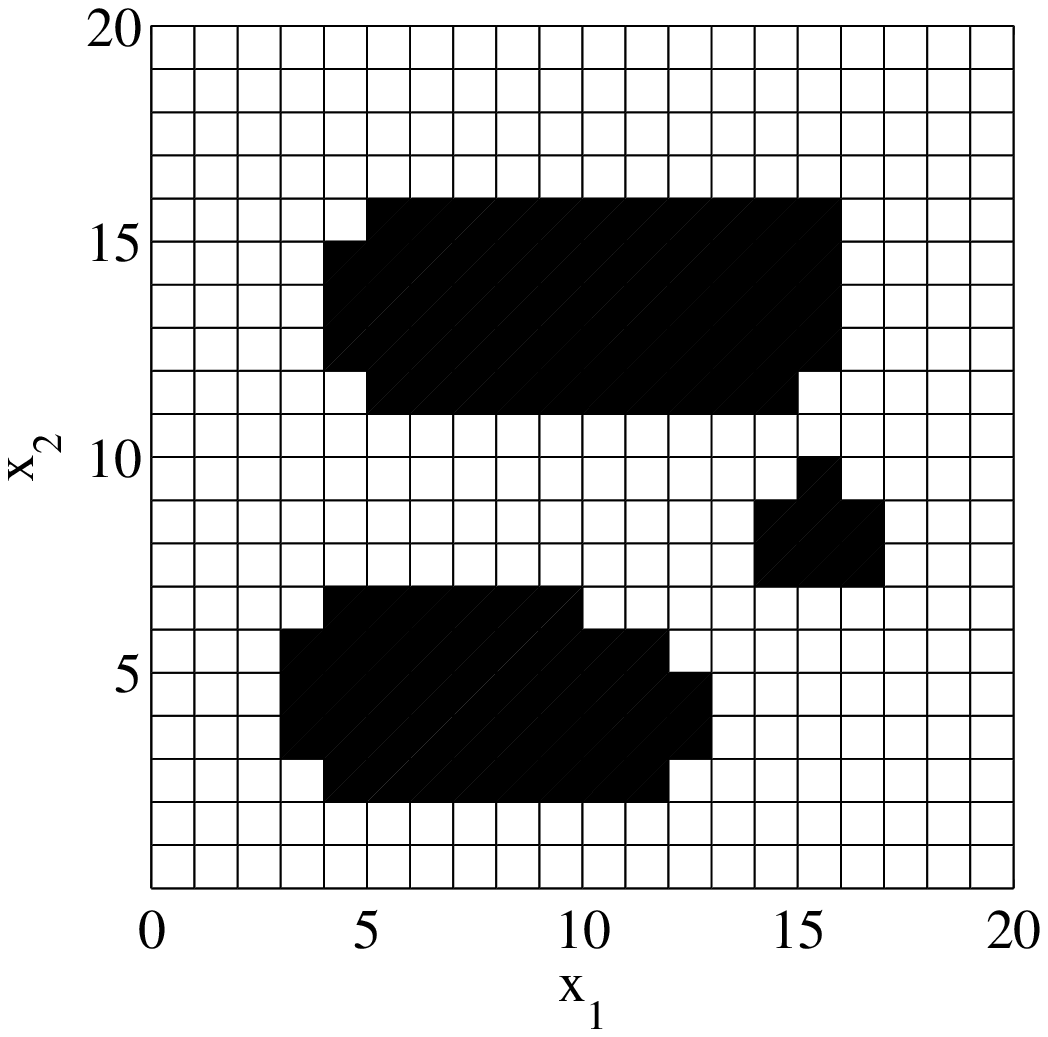}}
\subfigure[]{%
\includegraphics[width=6cm]{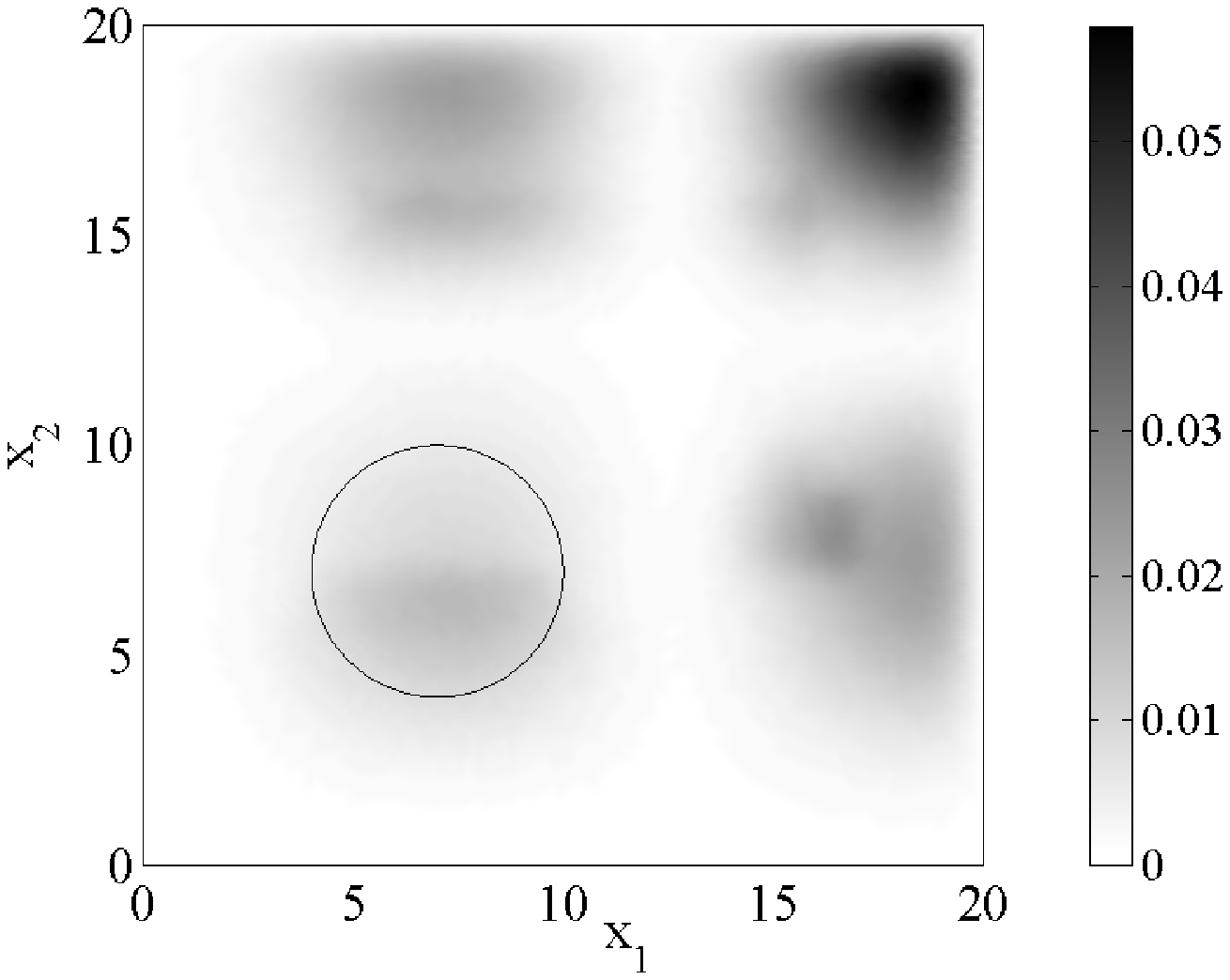}}
%\subfigure[]{%
%\includegraphics*[width=6cm]{exple3_mu2000r.eps}}
%\subfigure[]{%
%\includegraphics*[width=6cm]{exple3_mur.eps}}
\caption{Example 3: a) exact habitat configuration $\mu$. In the
black regions, the depicted function takes the value $2$; in the
white regions, it takes the value $-1$. b) measurement of $
u_{\gamma}(0.25,x)$ which was used for recovering $\mu$;  the
domain $\omega$ is delimited by a black circle. The exact
configuration $\mu$ was recovered after $3000$ iterations.
}%
\label{fig:exple3}
\end{figure}

We fixed $\mu\in E$ as in Fig. \ref{fig:exple3} (a). The measurement $ u_{\gamma}(\frac{t_0+t_1}{2},x)$ is depicted in Fig. \ref{fig:exple3} (b).
Two elements of $E$ are said to be neighbours  if they differ only on one cell of the lattice.

 The sequence
($\hat{\mu}_n$) stabilised on the exact configuration $\mu$ after
3000 iterations.

{\bf Example 4:} $E$ is composed of functions which can  take $21$ different  values between $m$ and $M$ on each cell of a regular lattice.

\begin{figure}
\centering
\subfigure[]{%
\includegraphics[width=6cm]{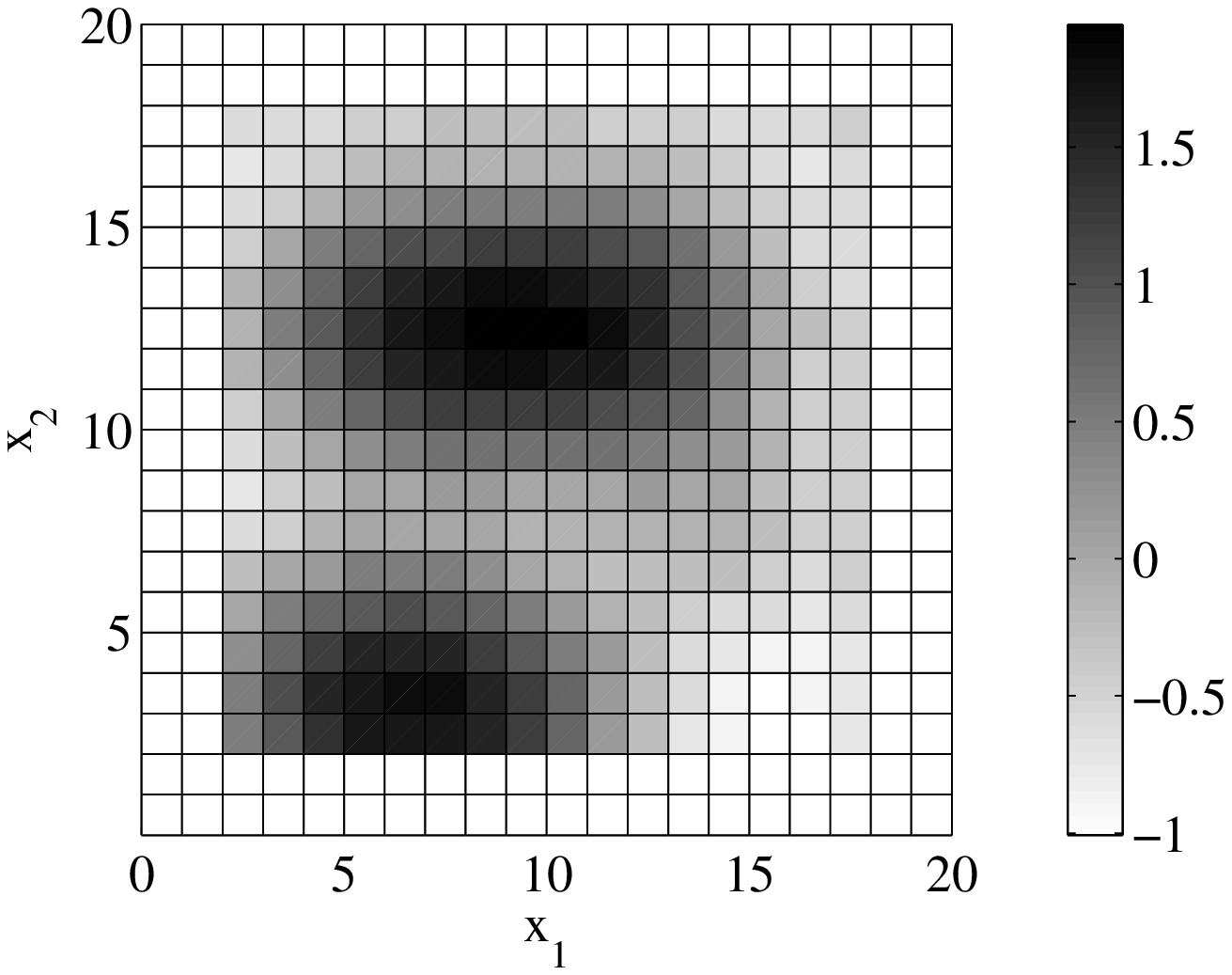}}
\subfigure[]{%
\includegraphics[width=6cm]{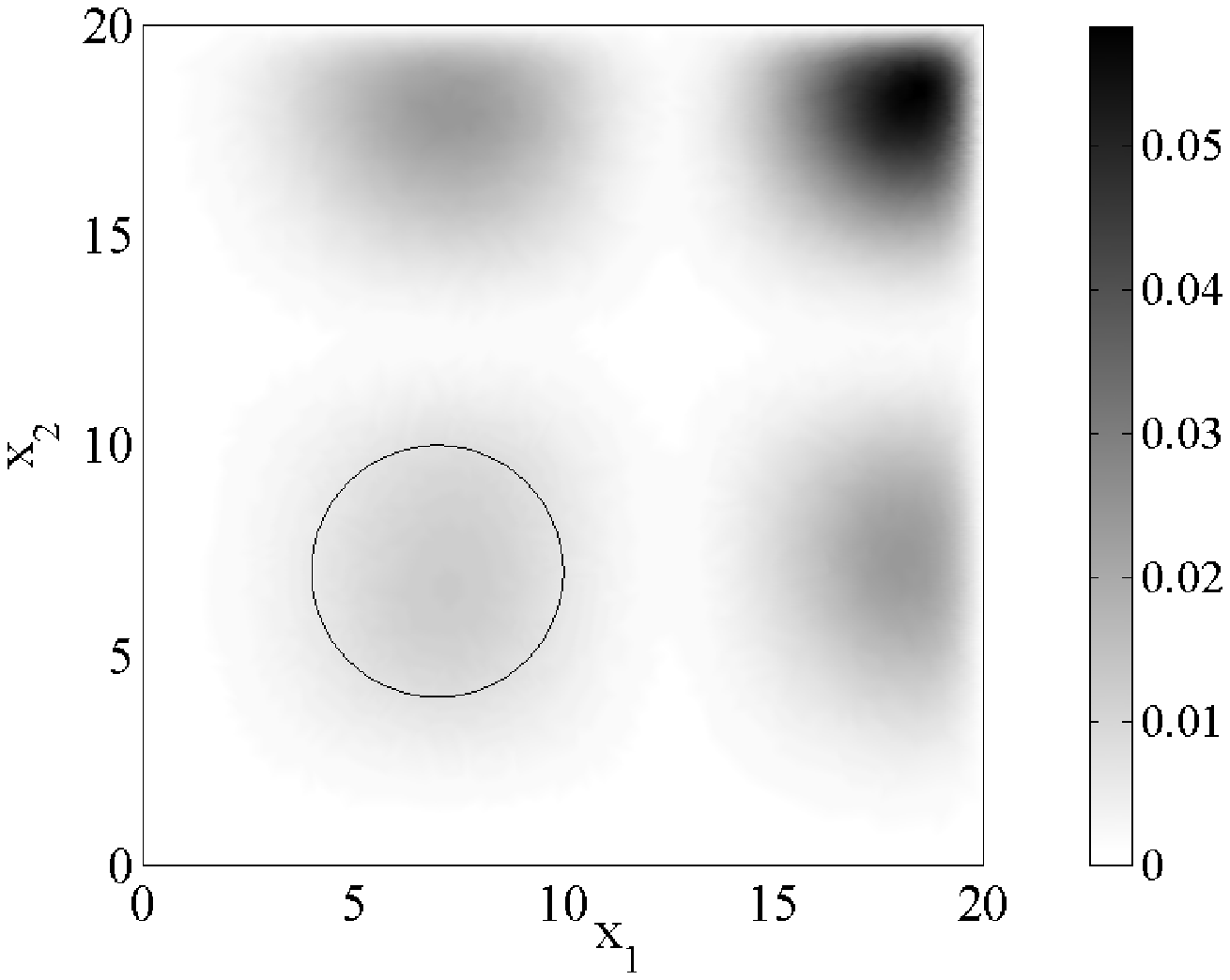}}
%\subfigure[]{%
%\includegraphics*[width=6cm]{exple4_mu7000r.eps}}
\subfigure[]{%
\includegraphics[width=6cm]{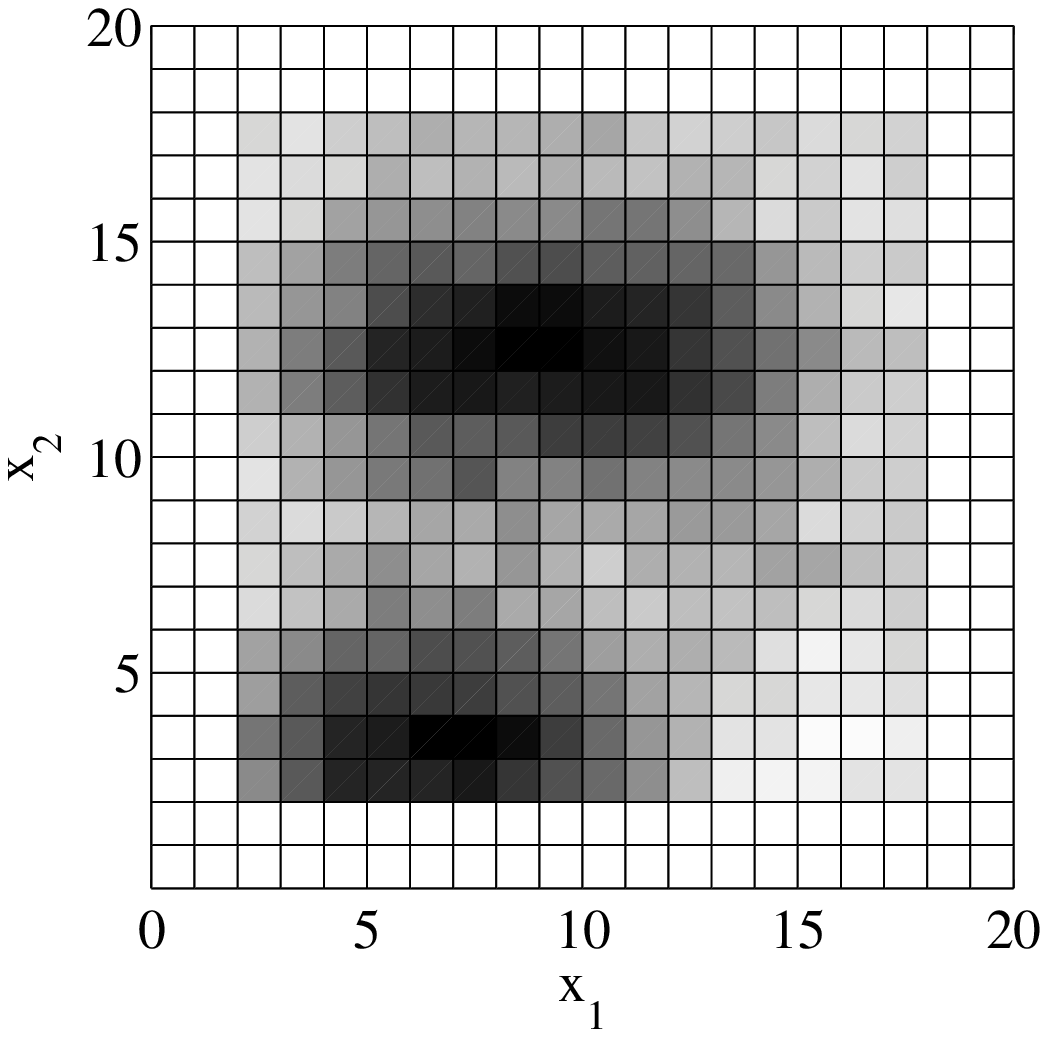}}
\caption{Example 4: a) exact habitat configuration $\mu$; b)
measurement of $ u_{\gamma}(0.25,x)$ which was used for recovering
$\mu$;  the domain $\omega$ is delimited by a black circle;  c)
approached configuration $\hat{\mu}_{14000}$.  The darker the
regions in (a), (b) and (c),  the higher the values of the
depicted functions.
}%
\label{fig:exple4}
\end{figure}

The configuration $\mu$ and the measurement $ u_{\gamma}(\frac{t_0+t_1}{2},x)$ are depicted in Fig. \ref{fig:exple4} (a) and (b).
Two elements of $E$ are said to be neighbours  if they differ by  $(M-m)/20$ on one cell of the lattice.

The sequence ($\hat{\mu}_n$)
stabilised on a configuration $\hat{\mu}$ after 14000 iterations
(Fig. \ref{fig:exple4}, (c)). The mean error  was
$\frac{1}{20^2}\int_0^{20}\int_0^{20}|\mu(x,y)-\hat{\mu}(x,y)|dx
dy=0.04.$

\section{Discussion and conclusion}\label{section5}

We have shown that, for an invasive species whose density is well
modelled by a reaction-diffusion equation, the spatial arrangement
of the favourable and unfavourable regions
can be measured indirectly through the
population density at the beginning of the invasion. More
precisely, we considered a logistic-like reaction-diffusion model,
and we placed ourselves under the assumption that the initial
population density was far from the environment carrying capacity
(it can be reasonably assumed at the beginning of an invasion). In
such a situation, the position of the favourable and unfavourable
regions, modelled through the intrinsic growth rate coefficient
$\mu$, may not be known \emph{a priori}. This is especially true for
exotic species whose hosts in native range and newly-colonised
areas can be different. From our results, in the ``ideal case"
considered here, the position of these regions can be obtained
through partial measurements of the population density. These
partial measurements consist in two samples of the population
density: 1) a ``spatio-temporal" measurement, but very locally (in
the small subset $\omega$) and during a short time period and, 2)
a ``spatial" measurement in the whole region susceptible to
colonisation ($\Omega$).

The stochastic algorithm presented in section \ref{section3} shows
explicitly how to reconstruct the habitat arrangement $\mu$ from
the above partial measurements of the population density. This
algorithm was proved to be effective in both one-dimensional and
two-dimensional cases, in section \ref{section4}, through several
numerical experiments. In examples 1 and 3, the algorithm converged
to the exact habitat configuration. In examples 2 and 4,
the sizes of the sets of possible habitat configurations were increased
compared to  examples 1 and 3. In those cases,  the algorithm
converged to configurations which were close to the exact ones.
It is noteworthy that the spatial measurement in $\Omega$ and the habitat arrangement $\mu$ can have very different shapes; therefore, $\mu$ cannot be straightforwardly deduced from this measurement.

These results can be helpful in preventing biological invasions.
Indeed, a
simple protocol, consisting of placing one trap in the invaded
region, and recording the number of individuals captured by this
trap over a short time-period (depending on the species
characteristics), and performing a single survey of the number of
individuals and their position in the whole considered region
should allow, from our results,  to detect the favourable areas,
and to treat them preventively. As we have emphasised in
Proposition \ref{lemp}.a, the knowledge of an estimate of the
habitat arrangement $\mu$ also allows us to forecast the final
population density, and therefore to detect the regions at higher
risk,  for instance in the case of harmful species.
As recalled in Proposition \ref{lemp}.b, having a good estimate of
the habitat arrangement $\mu$ is also crucial to forecast the
fate of the invasive species: persistence or extinction.

On the other hand, we have to underline that our approach may not be adapted to some species, especially those which colonies are made of few individuals. Indeed, the diffusion operator of our model can be obtained as the macroscopic limit of uncorrelated random walks. With such an operator, by the parabolic maximum principle \cite{evans}, it is known that, even with a compactly supported initial population density, the solution of our model is strictly positive everywhere on the domain as soon as $t>0$. This means that the solution, and therefore the information,  propagate  with infinite speed, which is not realistic for discrete populations. This could induce a practical limitation of our method to a certain type of species, which are well modelled by continuous diffusion processes  even at low densities (typically some insect or plant species, with high carrying capacity and growth rate).

Note also that some of the mathematical tools used in this paper, and especially
Carleman estimates (see Appendix A), were initially not adapted to
the nonlinear model considered here. Thus, we first considered, in
Theorem \ref{th1}, the linear - or Malthusian - case. For
populations whose density is far from the environment
carrying capacity, the linear and the nonlinear problem have close
solutions. In this situation, Theorem \ref{th2} extended the
result of Theorem \ref{th1} to the nonlinear case of a logistic
growth.

The results of this paper could be immediately extended to the
case of spatially varying functions $\gamma(x)$. Another easy
extension would be, for the ``spatio-temporal" measurement, to use
a partial boundary observation on a part $\Gamma^+$ of the domain
boundary $\partial\Omega$ instead of sampling the population over
a small domain $\omega$.   Using a new Carleman estimate (see
\cite{AlTa}) we are indeed able to write a stability result for
the coefficient $\mu$, similar to that of Theorem \ref{th1}, but
with $\|
\partial_{\mathbf{n}}(\partial_t u_{\gamma}) -
\partial_{\mathbf{n}}(\partial_t\tilde{v})
\|^2_{L^2((t_0,t_1)\times\Gamma^+)}$ instead of $\|\partial_t
u_{\gamma}-\partial_t\tilde{v}\|^2_{L^2((t_0,t_1)\times\omega)}$
in the definition of the functional $G$.

\section*{Acknowledgements} The authors would like to
thank the anonymous reviewers for their helpful comments and
insightful suggestions. This study was partly supported by the
french ``Agence Nationale de la Recherche" within the project
URTICLIM ``Anticipation des effets du changement climatique sur
l'impact \'ecologique et sanitaire d'insectes forestiers
urticants" and by the European Union within the FP 6 Integrated
Project ALARM- Assessing LArge-scale environmental Risks for
biodiversity with tested Methods (GOCE-CT-2003-506675).

\section{Appendices}

Let us introduce the following notations: for all $t, t'\in \R,$
with $t'>t$, we denote $Q_{t}^{t'}= (t,t')\times \Omega$ and
$\Sigma_{t}^{t'}= (t,t')\times
\partial \Omega$. Throughout this section, with a slight abuse of notation, we designate by
$C$ any upper bounds in our computations, provided they only
depend on the parameters $\Omega,$ $\Omega_1,$ $\omega,$
$\mathcal{B}_{\epsilon},$ $D,$ $t_0,$ $t_1$,
$\underline{u_i}/\overline{u_i}$.

\subsection{ Appendix A: proof of Theorem \ref{th1}}

{\bf Carleman estimate}

We recall here a Carleman-type estimate with a single observation.
Let $\beta$ be a function in
${\cal C}^2(\overline{\Omega})$ such that
$$
1<\beta<2 \mbox{ in }\Omega, \ \ \beta=1 \mbox{ on }\partial
\Omega, \ \min\{|\nabla  \beta (x)|, x \in
\Omega\backslash \overline{\omega}\}>0 \mbox{ and }{\partial}_{\bf
n}  \beta < 0 \mbox{ on }
\partial \Omega,
$$
where ${\bf n}$ denotes the outward unit normal to $\partial
\Omega$. For $\lambda> 0$ and $t \in (t_0,t_1)$, we define the
following weight functions
$$
 \label{wf}
 \varphi(t,x)=\frac{e^{\lambda \beta(x)}}{(t-t_0)(t_1-t)},
  \ \  \eta(t,x)=\frac{e^{2 \lambda } -e^{\lambda
 \beta(x)}}{(t-t_0)(t_1-t)}.
$$
Let $q$ be a solution of the parabolic problem
\begin{equation*} \left\{ \baa{l} \partial_t
q - D \Delta q + \alpha(x) q= f(t,x) \hbox{ in }Q_{t_0}^{t_1},\\
q =0 \hbox{ on }\Sigma_0^{t_1}, \\ q(0,x)=q_{0}(x) \hbox{ in
}\Omega, \eaa \right. \ (P)\end{equation*} for some functions
$f\in L^2(Q_{t_0}^{t_1}),$ and $\alpha, q_0\in L^{\infty}(\Omega)$.

Then the following results are proved in \cite{F:00}:

\begin{lemma}
\label{pr-Carl-Fur} Let $q$ be a solution of (P). Then, there
exist three positive constants $\lambda_0$, $C_0$ and $s>1$,
depending only on $\Omega$, $\omega$, $t_0$ and $t_1$ such that,
for any $\lambda \ge \lambda_0$, the next inequalities hold:
\begin{equation*}  \baa{l}  \hbox{a) }\|M_1(e^{-s
\eta}q)\|^2_{L^2(Q_{t_1})} + \|M_2(e^{-s
\eta}q)\|^2_{L^2(Q_{t_1})}+
s \lambda^2   \int_{Q_{t_0}^{t_1}} e^{-2s
\eta} \varphi |\nabla q|^2  \\  +s^3 \lambda^4
\int_{Q_{t_0}^{t_1}} e^{-2s \eta} \varphi^{3} q^2
 \leq C_0  \left[
  s^3 \lambda^{4}\int_{t_0}^{t_1}
\int_\omega e^{-2s \eta} \varphi^{3} q^2
 +   \int_{Q_{t_0}^{t_1}} e^{-2s \eta}\ ( f -\alpha q)^2 \right],
 \eaa
\end{equation*} where   $M_1$ and $M_2$ are defined by
$M_1\psi=-D \Delta\psi - s^{2}\lambda^{2} D
|\nabla\beta|^2\varphi^{2}\psi + s(\partial_t{\eta})\psi$, and
 $M_2\psi = \partial_t \psi +2 s \lambda D \varphi \nabla \beta . \nabla
\psi +2 s \lambda^{2} D \varphi | \nabla \beta |^2 \psi.$
Moreover,  \begin{equation*} \baa{l}
   \hbox{b) } s^{-2}\lambda^{-2}   \int_{Q_{t_0}^{t_1}} e^{-2s \eta} \varphi^{-1}
[(\partial_t q)^2+(\Delta q)^2]   +    \int_{Q_{t_0}^{t_1}} e^{-2s
\eta} \varphi |\nabla q|^2
 +s^2 \lambda^2   \int_{Q_{t_0}^{t_1}} e^{-2s \eta} \varphi^{3} q^2    \\
 \leq C_0  \left[
  s^2 \lambda^{2}\int_{t_0}^{t_1}
\int_\omega e^{-2s \eta} \varphi^{3} q^2
 + s^{-1} \lambda^{-2}    \int_{Q_{t_0}^{t_1}} e^{-2s \eta}\ ( f -\alpha q)^2
 \right].
\eaa \end{equation*}
\end{lemma}

{\bf Stability estimate with one observation}

Let $\mu, \tilde{\mu} \in \mathcal{M}$. We consider the solutions
$v$ and $\widetilde{v}$ of the linear problems $(P_{\mu,0})$ and
$(P_{\tilde{\mu},0}),$ respectively. We set $w=v-\widetilde{v}$,
$y=\partial_t w$,  and $\sigma= \mu-\widetilde{\mu}$. The function
$y$ is a solution of:
\begin{equation}
\left \{ \begin{array}{lll}
 \label{syst3}
   {\partial}_t y= D \Delta y+\mu y+\sigma \partial_t \widetilde{v} & \mbox{in} & Q_{t_0}^{t_1},\\
   y(t,x)=0& \mbox{on} & \Sigma_{t_0}^{t_1},\\
   y(0,x)=\sigma u_i(x)  & \mbox{in} & \Omega,\\
 \end{array}\right.
\end{equation}

The  function $\eta(x,t)$ attains  its
minimum value with respect to the time at $t= T'= \frac{t_0+t_1}{2}$. We set $\psi = e^{-s
\eta} y$. Using the operator $M_2$, introduced in Lemma
\ref{pr-Carl-Fur}, we introduce, following \cite{BP:02} (see also
\cite{CGR} and \cite{CCG}),
\begin{eqnarray*}
 \mathcal{I} =  \int_{t_0}^{ T' }\hspace*{-9pt}
 \int_{\Omega} M_2 \psi\;\psi.
\end{eqnarray*}
Let $\lambda_0$  be fixed as in Theorem
\ref{pr-Carl-Fur}.
\begin{lemma}
\label{lemma1} Let $\lambda \geq \lambda_0$. There exists a
constant $C$ such that
\begin{equation}
\label{eqlem1}
   |\mathcal{I}| \leq C
   \left[s^{3/2} \lambda^{2}  \int_{t_0}^{t_1}
 \int_{\omega} e^{-2s \eta} \varphi^{3}
    y^2
     + s^{-3/2} \lambda^{-2}
 \int_{Q_{t_0}^{t_1}}
   e^{-2s \eta}  \sigma^2 (\partial_t \widetilde{v})^2  \right].
 \end{equation}
\end{lemma}
\textit{Proof: } From the H\"older inequality, we have:
 \begin{equation*}
 |\mathcal{I}| \leq   s^{-3/2} \lambda^{-2}
\left( \int_{0}^{ T' } \int_{\Omega}
 (M_2 \psi)^2   \right)^{1/2}
\left( s^{3} \lambda^{4} \int_{0}^{ T' }  \int_{\Omega}
e^{-2s \eta}  y^2   \right)^{1/2}.
\end{equation*}
Thus using Young's inequality, we obtain
 \begin{equation}
   |\mathcal{I}| \leq   \frac{1}{4} s^{-3/2} \lambda^{-2}
   \left(\|M_2 \psi\|_{L^2(Q_{t_0}^{t_1})}^2 +
   s^3 \lambda^{4} \int_{Q_{t_0}^{t_1}}
    e^{-2s\eta} \varphi^3 y^2 \right).
    \label{young}
 \end{equation}
Applying  inequality a) of Lemma \ref{pr-Carl-Fur} to $q:=y$, we
obtain that there exists $C>0$, such that
 \begin{equation}
 \baa{l}
\|M_2 \psi\|_{L^2(Q_{t_0}^{t_1})}^2  +
  2 s^3 \lambda^{4}   \int_{Q_{t_0}^{t_1}}
    e^{-2s\eta} \varphi^3 y^2  \\  \leq C  \left[
  s^{3} \lambda^{4}\int_{t_0}^{t_1}
\int_\omega e^{-2s \eta} \varphi^{3} y^2
 +  \int_{Q_{t_0}^{t_1}} e^{-2s \eta}\ 2 [\mu^2 y^2 + \sigma^2 (\partial_t \widetilde{v})^2] \right]. \eaa
    \label{eqlem1_1}
 \end{equation}
Furthermore, since $\mu$ is bounded, and since $\varphi$ is
bounded from below by a positive constant, independent of
$\lambda$, we get that \be    \int_{Q_{t_0}^{t_1}} 2 e^{-2s \eta}
\mu^2 y^2   \leq s^3 \lambda^{4}  \int_{Q_{t_0}^{t_1}} e^{-2s\eta}
\varphi^3 y^2  , \label{eqlem1_2}  \ee for $\lambda$ large enough.
Combining (\ref{eqlem1_1}) and (\ref{eqlem1_2}), we obtain:
 \begin{equation}
 \baa{l}
\|M_2 \psi\|_{L^2(Q_{t_0}^{t_1})}^2  +
   s^3 \lambda^{4}   \int_{Q_{t_0}^{t_1}}
    e^{-2s\eta} \varphi^3 y^2    \\  \leq C  \left[
  s^{3} \lambda^{4}\int_{t_0}^{t_1}
\int_\omega e^{-2s \eta} \varphi^{3} y^2
 +\int  \int_{Q_{t_0}^{t_1}} e^{-2s \eta}\ \sigma^2 (\partial_t \widetilde{v})^2 \right]. \eaa
    \label{eqlem1_3}
 \end{equation}
 The conclusion of Lemma \ref{lemma1} follows from (\ref{young})
 and (\ref{eqlem1_3}). $\Box$

%%%%%%%%%%%%%%%%%%%%%%%%
% lemma                %
%%%%%%%%%%%%%%%%%%%%%%%%
\begin{lemma}
 \label{lemma2}
 Let $\lambda \geq \lambda_0$.
There exists a constant $C$ such that \begin{equation*}
 \begin{array}{ll}
 \label{eq:lemma2}
   \int_{\Omega}
   e^{{-2s \eta}( T' ,x)}\
   (\sigma \widetilde{v}( T' ,x))^2
  &  \leq C \left[ s^{3/2} \lambda^{2}  \int_{t_0}^{t_1}
 \int_{\omega} e^{-2s \eta} \varphi^{3}
    y^2  \right. \\
  & \left.  +   s^{-1} \lambda^{-2}\int_{Q_{t_0}^{t_1 }}
   e^{-2s \eta( T' ,x)}  \sigma^2 (\partial_t \widetilde{v})^2 \right. \\
   &\left. +   \int_{\Omega}
   e^{{-2s \eta}( T' ,x)}\
   (D \Delta\ w( T' ,x) +\mu w( T' ,x))^2 dx \right].
\end{array}
\end{equation*}
\end{lemma}
\textit{Proof: }  Using integration by parts over $\Omega$
and
the boundary condition $\psi=0$ on $\Sigma_{t_0}^{t_1}$, we
get:
\begin{equation}
\baa{ll}
 \mathcal{I}
 = & \frac{1}{2} \int_{Q_{t_0}^{ T' }} \partial_t (\psi^2)
\\ &
 - s \lambda D \int_{Q_{t_0}^{ T' }}
 \nabla  \cdot ({\varphi}
 \nabla \beta) \psi^2
 + 2s \lambda^2 D \int_{Q_{t_0}^{ T' }}
 {\varphi}
  |\nabla \beta|^2 \psi^2  . \eaa
\end{equation}
We then obtain
\begin{equation*}
 \label{eq:I}
 \frac{1}{2} \int_{\Omega} \psi\left( T' , \cdot\right)^2   = \mathcal{I}
 - s \lambda^2 D \int_{ Q_{t_0}^{T'}}
 {\varphi}
 \vert \nabla \beta \vert^2 \psi^2 dx
 dt + s \lambda D \int_{ Q_{t_0}^{T'}}{\varphi}\Delta \beta \psi^2\; dx dt
\end{equation*}
since $\psi(t_0)=0$ and $\nabla\varphi=\lambda\varphi\nabla\beta$.
As a consequence, for $\lambda > 1$,  and since $\Delta \beta$ is
bounded in $\Omega$, we finally get
\begin{equation} \label{eq:I2} \int_{\Omega}
 e^{{-2s \eta}( T' ,x)} y( T' ,x)^2 \; dx
\leq  2 |\mathcal{I}| + C s \lambda^2 \int_{t_0}^{ T' }
\int_{\Omega} e^{{-2s \eta}}{\varphi} y^2,
\end{equation}
for some constant $C$. Using \, $\varphi \leq
\frac{(t_1-t_0)^2}{4} \varphi^3$ and Lemma \ref{pr-Carl-Fur}, we
get that: \be s \lambda^2 \int_{Q_{t_0}^{T'}} e^{{-2s
\eta}}{\varphi} y^2 \leq C \left[
  \lambda^2 \int_{t_0}^{t_1}
\int_\omega e^{-2s \eta} \varphi^{3} y^2
 + s^{-2} \lambda^{-2}  \int_{Q_{t_0}^{t_1}} e^{-2s \eta}\ 2[\mu^2 y^2+ \sigma^2 (\partial_t \widetilde{v})^2]\right]. \label{eq:I2a}\ee
Arguing as in the proof of Lemma \ref{lemma1} for equation
(\ref{eqlem1_2}), and since, for all $x\in \Omega$, the function
$t\mapsto \eta(t,x)$ attains its minimum over $(t_0,t_1)$ at
$t=T'$, we finally obtain that, for $\lambda$ large enough, the
last term in (\ref{eq:I2}) is bounded from above by:
$$
C \left[  \lambda^{2}  \int_{t_0}^{t_1}
 \int_{\omega} e^{-2s \eta} \varphi^{3}
    y^2
    +  s^{-2} \lambda^{-2}\int_{Q_{t_0}^{t_1}}
   e^{-2s \eta( T' ,x)}  \sigma^2 (\partial_t
   \widetilde{v})^2\right].
$$
If we now observe that
$$y\left( T' ,x\right)= D \Delta w\left( T' ,x\right)+\mu w\left( T' ,x\right) + \sigma \partial_t\widetilde{v}\left( T' ,x\right),$$
we get:
$$  \sigma^2 \widetilde{v}\left( T' ,x\right)^2 \le   2 y\left( T' ,x\right)^2 + 2 \left[D \Delta w\left( T' ,x\right)+\mu w\left( T' ,x\right)\right]^2,$$
and, since $s>1$, the estimate of lemma \ref{lemma2} follows.
$\Box$

\begin{lemma}
We have
$0\leq v, \tilde{v} \leq \overline{u_{i}} e^{M t_1}$, and
$|\partial_t \tilde{v} | \leq (D+M)\overline{u_i}e^{M t_1}$ in
$Q_{t_0}^{t_1}$. \label{lem_estim_vt}
\end{lemma}
\textit{Proof: }From the parabolic maximum principle, we know that
$v, \tilde{v}\ge 0$ in $Q_{t_0}^{t_1}$.  Let $h$ be the solution of
the ordinary differential equation
\begin{equation} \left\{
\baa{l} h'=h M \hbox{ on }\R_+,\\
h(0)=\overline{u_{i}}. \eaa \right.\end{equation} The function $h$
is increasing and $H(t,\cdot):=h(t)$ is a supersolution of the
equations satisfied by $v$ and $\tilde{v}$. As a consequence of
the parabolic maximum principle, we have, \be 0 \leq v,\tilde{v}
\leq h(t_1)=\overline{u_{i}} e^{M t_1},\hbox{ in
}Q_{t_0}^{t_1}.\label{inegvv}\ee

Let us set $\rho:=\partial_t \tilde{v}$. The function $\rho$
satisfies:
\begin{equation}
\left \{ \begin{array}{lll}
 \label{systvt}
   {\partial}_t \rho= D \Delta \rho+\tilde{\mu} \rho & \mbox{in} & Q,\\
   \rho(t,x)=0& \mbox{on} & \Sigma,\\
   \rho(0,x)=D \Delta u_i(x) +\tilde{\mu} u_i  & \mbox{in} & \Omega.\\
 \end{array}\right.
\end{equation}
Since $u_i\in \mathcal{D}$, $\rho(0,x)\in L^{\infty}(\Omega)$.
Moreover, $H_1(t,x):=-(D+M) \overline{u_i}e^{M t}$ and
$H_2(t,x):=(D+M) \overline{u_i}e^{M t}$ are respectively sub- and
supersolutions of (\ref{systvt}). The parabolic maximum principle
leads to the inequalities,  \be -(D+M)\overline{u_i}e^{M t_1}\leq
\rho \leq (D+M)\overline{u_i}e^{M t_1} \hbox{ in
}Q_{t_0}^{t_1}.\label{inegvt}\ee
$\Box$

From Lemmas \ref{lemma2} and \ref{lem_estim_vt}, and since
$\sigma$ vanishes outside $\Omega_1$, it follows that
\begin{equation}
\begin{array}{l} \int_{\Omega_1} e^{-2s \eta( T' ,x)} \sigma^2
\left( \widetilde{v}(x, T' )^2 -  C  s^{-1} \lambda^{-2} (D+M)^2
\overline{u_i}^2 e^{2M t_1}
 \right) dx
\leq \\[2mm]
C s^{3/2} \lambda^{2} \int_{t_0}^{t_1}  \int_{\omega} e^{-2s \eta}
\varphi^{3}  y^2  + C \int_{\Omega} e^{{-2s \eta}( T' ,x)}
[(\Delta w( T' ,x))^2 + w( T' ,x)^2 ] \;dx.
\end{array} \label{eq:1}
 \end{equation}

\begin{lemma}
The ratio $\underline{u_i}/\overline{u_i}$ being fixed, there
exists $r>0$, independent of $\tilde{\mu}$, $u_i$,
$\underline{u_i}$ and $\overline{u_i}$, such that $\tilde{v}\left(
T' ,\cdot \right)>\overline{u_i} r$ in
$\Omega_1$.\label{lem_estiminf_v}
\end{lemma}
\textit{Proof: } Let $\xi_0 \leq 1$ be a smooth function in $\Omega$,
such that $\xi_0\equiv 1$ in $\mathcal{B}_{\frac{\epsilon}{2}}$
and  $\xi_0\equiv 0$ in
$\Omega\backslash\mathcal{B}_{\epsilon}$. Let $\xi$ be
the solution of
\begin{equation}
\left\{ \baa{l} \partial_t
\xi = D \Delta \xi -M \xi \hbox{ in }Q,\\
\xi(t,x)=0 \hbox{ on }\Sigma, \\
\xi(0,x)=(\underline{u_i}/\overline{u_i}) \xi_0(x) \hbox{ in
}\Omega. \eaa \right.
\end{equation}
Let us set $r:=\ds{\inf_{x\in \Omega_1}\xi\left( T' ,x\right)}$.
From the strong parabolic maximum principle, $\xi>0$ in
$Q_{t_0}^{t_1}$, and therefore, we get that $r>0$, since $\Omega_1$ is a closed subset of $\Omega$. Moreover, the parabolic maximum principle
also yields $\overline{u_i} \xi\leq \tilde{v}$ in $Q_{t_0}^{t_1}$.
In particular, we get $\tilde{v}\left( T' ,x\right)\geq
\overline{u_i} r>0$ in $\Omega_1$. $\Box$

From Lemma \ref{lem_estiminf_v}, it follows that, for
$\ds{\lambda\geq \sqrt{\frac{2C}{s}}\frac{(D+M) e^{M t_1}}{r}}$,
the term $\widetilde{v}(x, T' )^2 -  C s^{-1} \lambda^{-2} (D+M)^2
\overline{u_i}^2 e^{2M t_1}$ in (\ref{eq:1}) satisfies:
\begin{equation*}
\widetilde{v}(x, T' )^2 -  C s^{-1} \lambda^{-2} (D+M)^2
\overline{u_i}^2 e^{2M t_1} \geq  \overline{u_i}^2 \frac{r^2}{2}>0
\hbox{ in } \Omega_1.
\end{equation*}

We deduce that, for $\lambda$ large enough, \be \baa{l}
\int_{\Omega_1} e^{-2s \eta( T' ,x)} \sigma^2 dx  \\
\leq \frac{C}{\overline{u_i}^2 } \left(s^{3/2} \lambda^{2}
\int_{t_0}^{t_1} \int_{\omega} e^{-2s \eta} \varphi^{3}  y^2 +
\int_{\Omega} e^{{-2s \eta}( T' ,x)} ((\Delta w( T' ,x))^2 + w( T'
,x)^2) dx\right). \eaa \ee
Using the fact that $e^{-2s \eta} \varphi^{3}$ remains bounded in
$Q_{t_0}^{t_1}$,  we finally obtain \be
 \label{eq:sta}
 \|\sigma\|^2_{L^2(\Omega_1)} \leq \frac{C}{\overline{u_i}^2 } \left( \int_{t_0}^{t_1}  \int_{\omega}  y^2
+  \int_{\Omega} ((\Delta w( T' ,x))^2 + w( T' ,x)^2 ) \;dx
\right). \ee Recalling that $w=v-\widetilde{v}$, $y=\partial_t
(v-\widetilde{v})$,    (\ref{eq:sta}) implies  the result of
Theorem \ref{th1}.

\subsection{Appendix B: Proof of Theorem \ref{th2}}

Let $u_{\gamma}$ be the solution of $(P_{\mu,\gamma})$, and let
$v$ be the solution of $(P_{\mu,0})$. Let us set
$z_{\gamma}:=v-u_{\gamma}$. The function $z_{\gamma}$ is a
solution of
\begin{equation} \left\{ \baa{l}
\partial_t z_{\gamma}- D \Delta z_{\gamma} =z_{\gamma} \mu(x)+ \gamma u_{\gamma}^2\hbox{ in }Q,\\
z_{\gamma}(t,x)=0 \hbox{ on }\Sigma, \\
z_{\gamma}(0,x)=0 \hbox{ in }\Omega. \eaa
\right.\label{eqw}\end{equation} It follows from the parabolic
maximum principle that $z_{\gamma}(t,x)>0$ in $Q_{t_0}^{t_1}$. Thus
$u_{\gamma}\leq v$ in $Q_{t_0}^{t_1}$. Using the result of Lemma
\ref{lem_estim_vt}, we thus obtain: \be u_{\gamma}(t,x) \leq
\overline{u_{i}} e^{M t_1}\hbox{ in }Q_{t_0}^{t_1}.\label{inegu}\ee

Let $k$ be the solution of
\begin{equation} \left\{
\baa{l} k'=k M + \gamma \overline{u_{i}}^2 e^{2 M t_1} \hbox{ on }\R_+,\\
k(0)=0. \eaa \right.\end{equation} Since $K(t,x):=k(t)$ is a
supersolution of (\ref{eqw}), we obtain that $z_{\gamma}(t,x)\leq
k(t)$ in $Q_{0}^{t_1}$. Thus, since $k$ is increasing, \be
z_{\gamma}(t) \leq k(t_1)=\frac{\gamma \overline{u_{i}}^2 e^{2 M
t_1}}{M}(e^{M t_1}-1)\hbox{ in }Q_{t_0}^{t_1}. \label{inegw}\ee

Standard parabolic estimates (see e.g. \cite{evans}) then imply,
using (\ref{inegvv}), (\ref{inegu}), (\ref{inegw}) and the
hypothesis $\|u_i\|_{C^2(\overline{\Omega})}\leq \overline{u_i}$,
that:
\begin{eqnarray} \|\partial_t v\|_{L^2(Q_{t_0}^{t_1})} & = & \mathcal{O}(\|\mu v\|_{L^2(Q_{0}^{t_1})}+\|u_i\|_{H^1_0(\Omega)})= \mathcal{O}(\overline{u_i}), \label{estimv1} \\   \label{estimw1}
\|\partial_t z_{\gamma}\|_{L^2(Q_{t_0}^{t_1})} & = &
\mathcal{O}(\| z_{\gamma} \mu+ \gamma u_{\gamma}^2
\|_{L^2(Q_{0}^{t_1})}) = \mathcal{O}(\overline{u_{i}}^2),\\
\label{estimu1} \|\partial_t u_{\gamma}\|_{L^2(Q_{t_0}^{t_1})} & =
& \mathcal{O}(\| u_\gamma( \mu- \gamma u_{\gamma})
\|_{L^2(Q_{ 0}^{t_1})}+\|u_i\|_{H^1_0(\Omega)}) \nonumber \\ & = & \mathcal{O}(\overline{u_{i}}), \\
\sup_{t_0\leq t\leq t_1}\| v(t,.)\|_{H^2(\Omega)} & = &
\mathcal{O}(\|\mu v\|_{H^1(0,t_1;L^2(\Omega))}+
\|u_i\|_{H^2(\Omega)})= \mathcal{O}(\overline{u_i})
\label{estimv2} \\ \sup_{t_0\leq t\leq t_1}\|
z_{\gamma}(t,.)\|_{H^2(\Omega)} & = & \mathcal{O}(\|z_{\gamma}
\mu+ \gamma u_{\gamma}^2\|_{H^1(0,t_1;L^2(\Omega))}), \nonumber \\
& = & \mathcal{O}(\overline{u_i}^2) \label{estimw2}.
\end{eqnarray}
Using
(\ref{estimw1}) and (\ref{estimw2}), we get:
\be \label{estimw3}
G_{\mu}(\gamma,\mu) = \mathcal{O}(
\overline{u_{i}}^4).\ee Moreover, since
$u_0=v$, we have \be \baa{rl} G_{\mu}(\gamma,\tilde{\mu}) &
=G_{\mu}(0,\tilde{\mu})+G_{\mu}(\gamma,\mu)+2   \langle
\partial_t z_\gamma,\partial_t v- \partial_t \tilde{v} \rangle_{L^2(Q_{t_0}^{t_1})} \\
& +2   \langle \Delta z_\gamma\left( T' , \cdot\right),\Delta
v\left( T' , \cdot\right)- \Delta
\tilde{v}\left( T' , \cdot\right) \rangle_{L^2(\Omega)}\\
&+2   \langle z_\gamma\left( T' , \cdot\right), v\left( T' ,
\cdot\right)- \tilde{v}\left( T' , \cdot\right)
\rangle_{L^2(\Omega)}.
\eaa \ee
Thus, using (\ref{estimv1}) and (\ref{estimv2}), which are true
for both $v$ and $\tilde{v}$, and (\ref{estimw1}),
(\ref{estimw2}), (\ref{estimw3}), together with Cauchy-Schwarz
inequality, we get: \be
|G_{\mu}(\gamma,\tilde{\mu})-G_{\mu}(0,\tilde{\mu})|=\mathcal{O}(
\overline{u_i}^3). \ee  $\Box$

\subsection{Appendix C: Definitions of the state spaces $E$  and of their neighbourhood systems.}

In examples $1$ and $2$ of section $4$, $E$ is defined by
$$E:=\left\{\rho \in \mathcal{M}, \ \rho(x)=\sum_{k=10}^{89} \alpha_k \chi_k(x)\hbox{ in } \Omega_1, \hbox{ and }  \rho(x)=m \hbox{ in } \Omega\backslash \Omega_1 \right\},$$
where  $\chi_k$ are the characteristic functions of the intervals
$(\frac{k}{100},\frac{k+1}{100})$, and $\alpha_k$ are real numbers
taken in finite subsets of $[m,M]$.

In example $1$, $\alpha_k \in \{-1,2\}$, and two distinct elements $\nu_1$, $\nu_2$ of $E$, with
$\nu_1=\sum_{k=10}^{89} \alpha_{1,k} \chi_k(x)$ and
$\nu_2=\sum_{k=10}^{89} \alpha_{2,k} \chi_k(x)$ in $\Omega_1$,
are defined as neighbours if and only if   there exists a unique integer $k_0$ in  $[10,89]$ such that $\alpha_{1,k_0}\neq \alpha_{2,k_0}$. Note that, in this case, the number of elements in $E$ is $2^{80}$.

In example $2$,
$\alpha_k \in \{ (M-m) j/20+m, \hbox{ for } j=0\ldots 20.$\}, and  the neighbourhood system is defined as follows: two distinct
elements $\nu_1$, $\nu_2$ of $E$, with $\nu_1=\sum_{k=10}^{89}
\alpha_{1,k} \chi_k(x)$ and $\nu_2=\sum_{k=10}^{89} \alpha_{2,k}
\chi_k(x)$ in $\Omega_1$, are neighbours if and only if (i) there
exists a unique integer $k_0$ in  $[10,89]$ such that
$\alpha_{1,k_0}\neq \alpha_{2,k_0}$, (ii) additionally,
$|\alpha_{1,k_0}-\alpha_{2,k_0}|=(M-m)/20$. Note that,
in such a situation, the number of elements in $E$ is  $21^{80}$.

\

In examples $3$ and $4$ of section $4$, $E$ is defined by
$$E:=\left\{\rho \in \mathcal{M}, \ \rho(x)=\sum_{i,j=2}^{17} \alpha_{i,j} \chi_{i,j}(x)\hbox{ in } \Omega_1, \hbox{ and }  \rho(x)=m \hbox{ in } \Omega\backslash \Omega_1 \right\},$$
where  $\chi_{i,j}$ are the characteristic functions of the square
cells $(i/20,i/20+1/20)\times (j/20,j/20+1/20)$,  and
$\alpha_{i,j}$ are real numbers taken in finite subsets of
$[m,M]$.

In example $3$, $\alpha_{i,j} \in \{-1,2\}$. In this case,
the number of elements in $E$ is  $2^{256}$. Two distinct elements $\nu_1$, $\nu_2$ of $E$, with
$\nu_1=\sum_{i,j=2}^{17} \alpha^1_{i,j} \chi_{i,j}(x)$ and
$\nu_2=\sum_{i,j=2}^{17} \alpha^2_{i,j} \chi_{i,j}(x)$ for $x\in
\Omega_1$, are defined as neighbours if and only if there exists a
unique couple $(i_0,j_0)$ of integers  comprised between $2$ and $17$ such that
$\alpha^1_{ i_0,j_0 }\neq \alpha^2_{i_0,j_0}$.

In example $4$
$\alpha_{i,j} \in   \{ (M-m) j/20+m, \hbox{ for } j=0\ldots 20.$\}. The number of elements in $E$ is  ${21}^{256}$. In this case, two distinct elements $\nu_1$, $\nu_2$ of $E$, with
$\nu_1=\sum_{i,j=2}^{17} \alpha^1_{i,j} \chi_{i,j}(x)$ and
$\nu_2=\sum_{i,j=2}^{17} \alpha^2_{i,j} \chi_{i,j}(x)$ for $x\in
\Omega_1$, are defined as neighbours if and only if (i)
it exists a unique couple $(i_0,j_0)$ of integers  comprised between $2$ and $17$ such that $\alpha^1_{ i_0,j_0
}\neq \alpha^2_{ i_0,j_0 }$; and (ii) additionally
$|\alpha^1_{i_0,j_0}-\alpha^2_{ i_0,j_0 }|=(M-m)/20$.

\end{document}